\documentclass[11pt]{amsart}

\usepackage{enumerate}
\usepackage{graphicx}
\usepackage{soul}
\usepackage[usenames, dvipsnames]{color}
\usepackage{fancyhdr}
\usepackage{float}
\usepackage{bbm}
\usepackage[margin=1in]{geometry}
\usepackage[english]{babel}
\usepackage[protrusion=true,expansion=true]{microtype}	
\usepackage{amsmath,amsfonts,amsthm,amssymb,hyperref}
\usepackage{graphicx}
\usepackage{enumitem}
\usepackage{xpatch}
\newcommand{\N}{{\mathbb N}}
\def \e {\epsilon}
\def \O {\Omega}
\def \l {\lambda}

\def \r {\mathcal R}
\def \R {\mathbb{R}}

\usepackage{color}
\usepackage{xcolor}
\newtheorem{Theorem}{Theorem}[section]

\newtheorem{Lemma}[Theorem]{Lemma}

\newtheorem{Proposition}[Theorem]{Proposition}

\newtheorem{Remark}[Theorem]{Remark}

\def \N {\mathbb{N}}
\def \R {\mathbb{R}}

\def \l {\lambda}
\def \l {\lambda}
\def \B {\mathcal{B}}
\def \GN {\Gamma_N}
\def \GD {\Gamma_D}
\newcommand{\norm}[1]{\left\Vert #1 \right\Vert}
\newcommand{\normltwo}[2][\O]{\norm{#2}_{L^2(#1)}}
\newcommand{\inprod}[2]{\left\langle #1, #2 \right\rangle}
\newcommand{\inprodb}[2]{\B\left[#1,#2\right]}%{\left\langle #1, #2 \right\rangle_{\B}}
\newcommand{\inprodltwo}[3][\O]{\inprod{#2}{#3}_{L^2(#1)}}
\newcommand{\normb}[1]{\norm{#1}_\B}

\newcommand{\nderiv}[1][\nu]{\partial_#1}

\definecolor{punk}{rgb}{0.870, .678, 0.8}

\newcommand{\normhone}[2][\O]{\norm{#2}_{H^1(#1)}}

\makeatletter
\newcounter{proofpart}
\xpretocmd{\proof}{\setcounter{proofpart}{0}}{}{}
\newcommand{\proofpart}{%
  \par
  \stepcounter{proofpart}%
  \noindent\emph{Part \theproofpart: }%\nobreak\smallskip
}
\makeatother

\title[]{Inverse Iteration for the Laplace Eigenvalue Problem With Robin and Mixed Boundary Conditions}
\author[B. Lyons]{Benjamin Lyons} 
\address{Department of Mathematics, Rose-Hulman Institute of Technology, 5500 Wabash Ave., Terre Haute, IN 47803}
\email{lyonsba1@rose-hulman.edu}

\author[E. Ruttenberg]{Emily Ruttenberg}
\address{Department of Mathematics, University of Maryland Baltimore County, 1000 Hilltop Cir., Baltimore, MD 21250}
\email{eruttenberg@umbc.edu}

\author[N. Zitzelberger]{Nicholas Zitzelberger}
\address{Department of Mathematics, Oregon State University, 1500 SW Jefferson Way, Corvallis, OR 97331}
\email{zitzelbn@oregonstate.edu}



\begin{document}

\begin{abstract}
    We apply the method of inverse iteration to the Laplace eigenvalue problem with Robin and mixed Dirichlet-Neumann boundary conditions, respectively. For each problem, we prove convergence of the iterates to a non-trivial principal eigenfunction and show that the corresponding Rayleigh quotients converge to the principal eigenvalue. We also propose a related iterative method for an eigenvalue problem arising from a model for optimal insulation and provide some partial results. 
\end{abstract}

\maketitle
\footnotetext{This article was published by the PUMP Journal of Undergraduate Research, viewable at: \url{https://journals.calstate.edu/pump/index}}

\section{Introduction}
\subsection{The Dirichlet Eigenvalue Problem}
Let $\O \subset \R^n$, $n \geq 2$, be a smooth and bounded domain. The \emph{Dirichlet eigenvalue problem} for the standard Laplace operator on $\O$ is to find a pair $(\lambda, u) \in \R \times H^1_0(\O)$ satisfying
\begin{align}\label{DirichletEigenvalueProblem}
    \begin{cases}
        -\Delta u = \lambda u &\text{ in } \O, \\
        u = 0 & \text{ on } \partial \O.
    \end{cases}
\end{align}
We recall some well-known facts about \eqref{DirichletEigenvalueProblem}, referring to \cite[Section 6.5]{evans} and \cite[Section 7.41]{medkova} for details. The eigenvalues of \eqref{DirichletEigenvalueProblem} are discrete and strictly positive \cite[Section 6.5]{evans}, and the principal (smallest, non-zero) eigenvalue, which we denote $\lambda_D^1$, admits a variational characterization in terms of the Rayleigh quotient
\begin{align}\label{MixedRayleighQuotient}
    \r(v) := \frac{\int_\O |\nabla v|^2 \, dx}{\int_\O v^2 \, dx}.
\end{align}
Specifically,
\[
    \lambda_D^1 = \min_{v \in H^1_0(\O)} \r(v).
\]
Any eigenfunction corresponding to the eigenvalue $\lambda_D^1$ must be either strictly positive or strictly negative in $\O$. Furthermore, $\lambda_D^1$  is simple, i.e. any two eigenfunctions with eigenvalue $\lambda_D^1$ are scalar multiples of each other.

\subsection{The Robin Eigenvalue Problem}\label{RobinIntro}
Given a function $h \in C^1(\partial\O)$, the \emph{Robin eigenvalue problem} is to find a pair $(\lambda, u) \in \R \times H^1(\O)$ satisfying
\begin{align}\label{RobinEigenvalueProblem}
    \begin{cases}
        -\Delta u = \lambda u & \quad \text{in } \O, \\
        u + h\partial_\nu u = 0 & \quad \text{on } \partial \O,
    \end{cases}
\end{align}
where $\nu$ is the outward-pointing normal vector to $\partial \O$ and $\partial_\nu u = \nabla u \cdot \nu$ is the normal derivative of $u$. Note that when $h \equiv 0$, we recover the Dirichlet eigenvalue problem \eqref{DirichletEigenvalueProblem}. We will assume, from here onward, that $h$ is strictly positive on $\partial \O$. As we discuss in Section \ref{RobinBackgroundSec}, the Robin eigenvalue problem \eqref{RobinEigenvalueProblem} also has a discrete spectrum of strictly positive eigenvalues. The principal eigenvalue, which we denote $\lambda_R^1$, has a variational characterization in terms of the Rayleigh quotient 
\begin{align}\label{RobinRayleighQuotient}
    R(v) := \frac{\int_\O |\nabla v|^2 \, dx + \int_{\partial \O} \frac{v^2}{h}\, d\sigma}{\int_\O v^2 \, dx}.
\end{align}
Specifically,
\begin{align*}
    \lambda_R^1 = \min\limits_{v \in H^1(\O)} R(v).
\end{align*}
Note that the boundary integral appearing in the Rayleigh quotient $R(v)$ is well-defined by the fact that $h$ is strictly positive on $\partial \O$ and by the trace theorem\footnote{In order to keep our notation concise, we will not explicitly denote the trace operator in boundary integrals throughout this work.} for Sobolev functions \cite[Section 6.5]{evans}.
Finally, similar to the principal eigenfunctions for the Dirichlet problem, any eigenfunction corresponding to the principal Robin eigenvalue $\lambda_R^1$ must be either strictly positive or strictly negative in $\O$, and $\lambda_R^1$ is simple, i.e. any two eigenfunctions with eigenvalue $\lambda_R^1$ are scalar multiples of each other.

\begin{figure}[H]
    \centering
\includegraphics[width=0.5\linewidth]{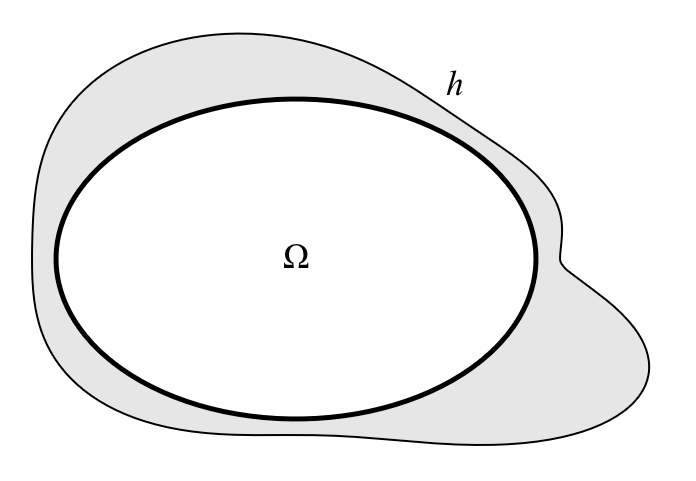}
    \caption{A conductor $\O$ with insulation $h$.}
    \label{fig:RobinFigure}
\end{figure}

The eigenvalue problem \eqref{RobinEigenvalueProblem} is tied to the heat diffusion problem on a homogeneous, conductive medium $\O$ surrounded by an insulating layer described by the function $h$ (see Figure \ref{fig:RobinFigure}). More specifically, let $v(x,t)$ be the solution of the heat equation with Robin boundary condition
\begin{align*}
    \begin{cases}
        v_t - \Delta v = 0 & \quad \text{in } \O, \\
        v + h\partial_\nu v = 0 & \quad \text{on } \partial \O, \\
        v(x,0) = v_0(x).
    \end{cases}
\end{align*}
The solution $v(x,t)$ represents the temperature distribution at $x \in \O$ at time $t$ given the insulation $h$ and initial temperature distribution $v_0(x)$. The long-time behavior of $v(x, t)$ is governed by the eigenvalues in \eqref{RobinEigenvalueProblem}; indeed, $v(x,t)$ approaches equilibrium with an asymptotic rate $e^{-\lambda_R^1 t}$ as $t \to \infty$.

\subsection{Mixed Boundary Conditions}\label{MixedIntro}

\begin{figure}[H]
    \centering
\includegraphics[width=0.4\linewidth]{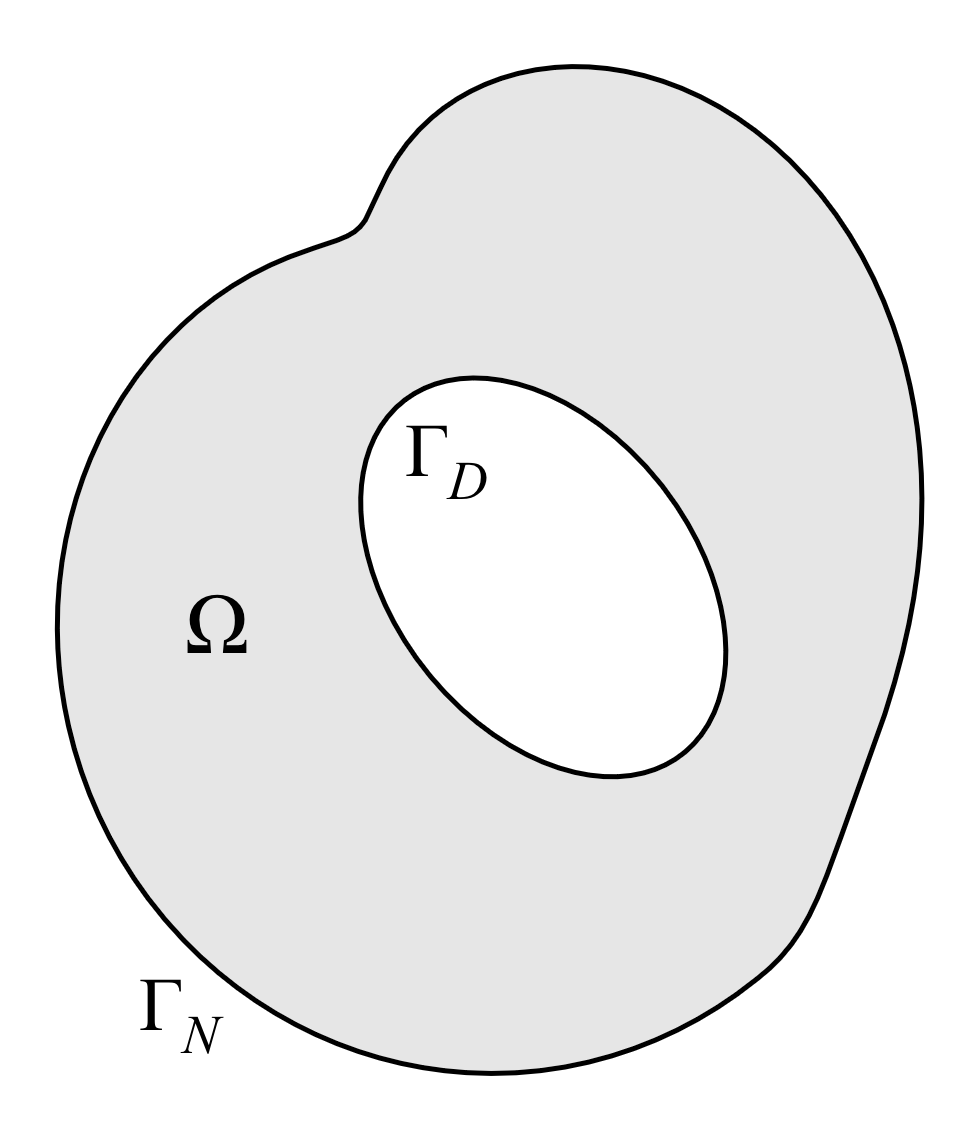}
    \caption{A doubly connected domain $\O$ with two connected boundary components $\GD$ and $\GN$.}
    \label{fig:doubledomain}
\end{figure}

Suppose now that $\partial\O$ consists of two disjoint, closed, connected components $\Gamma_D$ and $\Gamma_N$ (i.e. $\partial \O = \Gamma_D \cup \Gamma_N$ and $\Gamma_D \cap \Gamma_N = \emptyset$, as in Figure \ref{fig:doubledomain}). Fix $\e > 0$ and consider the Robin eigenvalue problem \eqref{RobinEigenvalueProblem} with piecewise constant insulation
\[h(x) = \begin{cases}
    \e &\text{ on } \Gamma_D,\\
    \frac{1}{\e} &\text{ on } \Gamma_N.
\end{cases}\] 
Taking a formal limit as $\e \to 0^+$, we obtain the \emph{mixed eigenvalue problem}, which is to find $(\lambda, u) \in \mathbb{R} \times \mathcal{C}$ satisfying
\begin{align}\label{MixedEigenvalueProblem}
    \begin{cases}
        -\Delta u = \lambda u & \quad \text{in } \O, \\
        u = 0 & \quad \text{on } \GD, \\
        \partial_\nu u = 0 & \quad \text{on } \GN,
    \end{cases}
\end{align}
where
\begin{align}\label{MixedTestFunctionClass}
    \mathcal{C} := \{v \in H^1(\O) : v\vert_{\Gamma_D} = 0\}.
\end{align}

As discussed in Section \ref{MixedBackgroundSec}, the eigenvalue problem \eqref{MixedEigenvalueProblem} also has a discrete spectrum of strictly positive eigenvalues. The principal eigenvalue, which we denote $\lambda_M^1$, has a variational characterization in terms of the Rayleigh quotient \eqref{MixedRayleighQuotient} given by
\begin{align*}
    \lambda_M^1 = \min\limits_{v \in \mathcal{C}}\r(v).
\end{align*}

Moreover, any eigenfunction corresponding to the eigenvalue $\lambda_M^1$ must be either strictly positive or strictly negative in $\O$, and $\lambda_M^1$ is simple, i.e. any two eigenfunctions with eigenvalue $\lambda_M^1$ are scalar multiples of each other.

\subsection{Main Results}
For physical applications, it is important to accurately approximate the principal eigenvalue and corresponding eigenfunctions for the problems \eqref{DirichletEigenvalueProblem}, \eqref{RobinEigenvalueProblem} and \eqref{MixedEigenvalueProblem}. For \eqref{DirichletEigenvalueProblem}, an approximation procedure that has been investigated by several authors (see Subsection \ref{litrev} below for specific references) is \emph{inverse iteration}. At a high level, this method generates a sequence of functions and positive real numbers converging to the principal eigenfunction and eigenvalue respectively by repeatedly updating and solving appropriate boundary-value problems for the Poisson equation. The appeal of such an iterative method is that  solving eigenvalue problems becomes equivalent to solving a (countably infinite) family of Poisson problems. Our main results concern the convergence of analogous inverse iteration schemes for computing the principal eigenvalues and eigenfunctions for both the Robin and mixed eigenvalue problems \eqref{RobinEigenvalueProblem} and \eqref{MixedEigenvalueProblem}. We follow the proof strategy outlined in \cite{MongeAmpereInverseIteration, NamLe} for a different class of eigenvalue problems.

\begin{Theorem}\label{RobinConvergence}
Let $u_0 \in C^1(\O) \cap L^\infty(\O)$ be a strictly positive function on $\overline{\O}$. For each $k \in \N$, let $u_{k+1}$ be the unique, classical solution of \begin{equation}\label{LaplacianIterationinMainResults}\begin{cases}
        -\Delta u_{k+1} = R(u_k) u_k & \quad \text{in } \O, \\
        u_{k+1} + h\partial_\nu u_{k+1}= 0 & \quad \text{on } \partial \O.
\end{cases}\end{equation}
Then $\lim\limits_{k \to \infty} R(u_k) = \lambda^1_R$ and $u_k$ converges strongly in $H^1(\O)$ to a positive principal eigenfunction for the Robin eigenvalue problem \eqref{RobinEigenvalueProblem}.
\end{Theorem}
\begin{Theorem}\label{MixedConvergence}
Let $u_0 \in C^1(\O) \cap L^\infty(\O)$ be a strictly positive function on $\overline{\O}$ satisfying $\r(u_0) > 0$. For each $k \in \N$, let $u_{k+1}$ be the unique, classical solution of \begin{equation}\label{MixedIteration}
        \begin{cases}
            -\Delta u_{k+1} = \r(u_k) u_k & \quad \text{in } \O, \\
            u_{k+1} = 0 & \quad \text{on } \GD, \\
            \partial_\nu u_{k+1} = 0 & \quad \text{on } \GN. 
        \end{cases}
    \end{equation}
Then $\lim\limits_{k \to\infty} \r(u_k) = \lambda^1_M$ and $u_k$ converges strongly in $H^1(\O)$ to a positive principal eigenfunction for the mixed eigenvalue problem \eqref{MixedEigenvalueProblem}.
\end{Theorem}

\begin{Remark}
Interestingly, the convergence of the Rayleigh quotients of the iterates to the principal eigenvalue can be shown \emph{before} establishing convergence of the sequence of iterates to an eigenfunction (see Propositions \ref{RayleighQuotientConvergence} and \ref{RayleighQuotientConvergenceMixed}). In fact, the convergence of the Rayleigh quotients plays an essential role in the proof of convergence of the iterates.
\end{Remark}

We also investigate an inverse iteration scheme related to an eigenvalue problem arising from the optimal insulation of conducting bodies studied in \cite{BucurButtazzo}.  While we obtain some positive results about our proposed iteration, we are unable to establish convergence to a solution of the optimization problem. We refer to Section \ref{sec:Double} for more details on the results we obtained and the difficulties we encountered.

\subsection{Comparison With Existing Literature}\label{litrev} As stated above, inverse iteration for the principal Dirichlet eigenvalue problem has been studied extensively \cite{MongeAmpereInverseIteration, Bozorgnia, NamLe, BiezunerDirichletIteration, HyndLindgrenProceedings}. In fact, each of these works considers the Dirichlet eigenvalue problem for various elliptic operators that generalize standard Laplacian. Inverse iteration for the Robin eigenvalue problem with constant $h$ has been considered in \cite{HyndLindgrenJFA, ConstantHPaper}. Inverse iteration for the mixed boundary-value problem has not, to the best of our knowledge, been studied explicitly in any prior work. 

Let us highlight some of the differences between the works \cite{HyndLindgrenJFA, ConstantHPaper} and the present manuscript. The authors of \cite{ConstantHPaper} consider the case where the insulation $h$ is identically equal to a positive constant (denoted $\beta^{-1}$ in \cite{ConstantHPaper}). Moreover, their sequence of iterates are normalized at each step to have unit $L^2$ norm. Our work  generalizes the results of \cite{ConstantHPaper} to the case of positive $h \in C^1(\partial \O)$, as well as the (formal) limiting case of the mixed eigenvalue problem. We also avoid any additional normalization step. Let us note that while \cite{ConstantHPaper} demonstrates the existence of a subsequence of iterates converging to a principal eigenfunction, the convergence of the full sequence of iterates is not properly addressed. We provide a complete proof of this convergence in Section \ref{RobinProofSec}. 

The work \cite{HyndLindgrenJFA} investigates an inverse iteration scheme for a class of abstract variational problems involving homogeneous functionals on Banach spaces that includes the Dirichlet and Robin eigenvalue problems as specific examples. For the Robin boundary condition (see \cite[Examples 2.11 and 3.12]{HyndLindgrenJFA}), the iteration considered in \cite{HyndLindgrenJFA} is
\begin{equation}\label{HyndLindgrenIteration}
    \begin{cases}
        -\Delta v_{k+1} = v_k &\text{ in } \O,\\
        v_{k+1} + h \nderiv v_{k+1} = 0 &\text{ on } \partial\O.
    \end{cases}
\end{equation}
Note the absence of a Rayleigh quotient term on the right-hand side, as compared to \eqref{LaplacianIterationinMainResults}. One of the main results of \cite{HyndLindgrenJFA} is that if the sequence of functions $(\lambda^1_R)^k v_k$ converges to a non-trivial function $v$, then $v$ is a principal eigenfunction for the Robin problem. Moreover, if $v$ is non-trivial, then
\[\lambda^1_R = \lim_{k \to \infty} \frac{\normltwo{v_{k-1}}}{\normltwo{v_k}}.\]
As discussed in \cite{HyndLindgrenJFA}, $(\lambda^1_R)^k v_k$ may converge to the zero function, depending on the choice of $v_0$; a sufficient condition is for $v_0$ to be bounded below by a positive principal eigenfunction (see discussion after displayed equation (2.4) in \cite{HyndLindgrenProceedings}). By comparison, our work shows that, under mild hypotheses on $u_0$, the sequence of inverse iterates  converge to a non-trivial principal eigenfunction, and the corresponding Rayleigh quotients converge to the principal eigenvalue. No prior knowledge of the eigenvalue or eigenfunctions is necessary to obtain a non-trivial limit.

Let us comment on how the sequence of iterates $\{u_k\}_{k=0}^\infty$ defined as in Theorem \eqref{RobinConvergence} are related to the sequence $\{v_k\}_{k=0}^\infty$ defined as in \eqref{HyndLindgrenIteration} and studied in \cite{HyndLindgrenJFA}. 
By unique solvability of the Poisson problem and an induction argument, it is not difficult to see that
\[(\lambda_R^1)^k v_k = \frac{u_k}{\prod\limits_{j=0}^{k - 1} \left(\frac{R(u_j)}{\lambda_R^1}\right)}.\]
Since Theorem \ref{RobinConvergence} shows $u_k$ converges to a non-trivial eigenfunction for a large family of initial choices $u_0$, the convergence of $(\lambda_R^1)^k v_k$ to a non-trivial function hinges on the convergence of the infinite product $\prod_{j=0}^{\infty} \frac{R(u_j)}{\lambda_R^1}$, which holds if and only if  
\begin{equation}\label{HyndLindgrenSum}    \sum_{j=0}^\infty\left( R(u_j) - \lambda_R^1\right) < \infty.
\end{equation}
While our work shows $R(u_j) - \lambda_R^1 \searrow 0$ as $j \to \infty$, we do not know how fast this convergence is occurring and  it is unclear what properties of $u_0$ will ensure \eqref{HyndLindgrenSum} holds. Indeed, establishing \emph{any} rate of convergence of the Rayleigh quotients to the principal eigenvalue appears to be a challenging problem.

\subsection{Outline}\label{Outline}

In Section \ref{sec:Robin}, we will study the Robin eigenvalue problem \eqref{RobinEigenvalueProblem}, establish some basic properties of the sequence generated by \eqref{LaplacianIterationinMainResults}, and prove Theorem \ref{RobinConvergence}. In Section \ref{sec:Mixed}, we will provide a similar analysis for the sequence generated by \eqref{MixedIteration} and prove Theorem \ref{MixedConvergence}. Finally, in Section \ref{sec:Double} we will investigate an iterative scheme relating to a nonlinear eigenvalue problem in optimal insulation and provide some partial results.

\section{Inverse Iteration for the Robin Eigenvalue Problem}\label{sec:Robin}
\subsection{Background}\label{RobinBackgroundSec}
By adapting the techniques in, for instance, \cite[Section 6.5, Theorems 1-2]{evans} and \cite[Theorems 7.41.5-7.41.6]{medkova}, one can prove the following:

\begin{enumerate}%[label = (\roman*)]
    \item The eigenvalues of \eqref{RobinEigenvalueProblem} are discrete, positive, and diverge to infinity.

    \item The principal eigenvalue $\lambda_R^1$ admits the variational characterization
    \begin{align*}
        \lambda_R^1 = \min\limits_{v \in H^1(\O)}R(v),
    \end{align*}
    where the Rayleigh quotient $R(v)$ is defined as in \eqref{RobinRayleighQuotient}.

    \item The principal eigenfunctions have a sign on $\O$.
    
    \item The principal eigenvalue $\lambda_R^1$ is simple; that is, any two principal eigenfunctions are scalar multiples of each other.

\end{enumerate}
If a function $u \in H^1(\O)$ satisfies 
% is a weak solution of \eqref{RobinEigenvalueProblem} with $\lambda = \lambda_R^1$ if
\begin{equation}\label{RobinEigenvalueWeakFormulation}
    \int_{\O}\nabla u \cdot \nabla \varphi \, dx + \int_{\partial\O} \frac{u\varphi}{h} \, d\sigma = \lambda_R^1\int_{\O}u\varphi \, dx \qquad \forall\varphi \in H^1(\O),
\end{equation}
then $u$ is called a weak solution of the Robin eigenvalue problem \eqref{RobinEigenvalueProblem} with $\lambda = \lambda_R^1$.

Consider the bilinear form on $H^1(\O)$ defined by
\begin{align*}
    \inprodb{u}{v} := \int_{\O}\nabla u \cdot \nabla v \, dx + \int_{\partial\O} \frac{uv}{h} \, d\sigma.
\end{align*}
Since $h$ is strictly positive and $\partial \O$ is compact, it follows from \cite[Theorem 7.43.1]{medkova} that $\inprodb{\cdot}{\cdot}$ is an inner product, and its induced norm $\normb{\cdot}$ is equivalent to the $H^1(\O)$ norm. We can recast \eqref{RobinEigenvalueWeakFormulation} as
\begin{align}\label{RobinEigenvalueWeakFormulationShorthand}
    \inprodb{u}{\varphi} = \lambda_R^1\inprodltwo{u}{\varphi} \qquad \forall\varphi \in H^1(\O)
\end{align}
and the Rayleigh quotient can be written as
\begin{align}\label{RobinRayleighQuotientShorthand}
    R(v) = \frac{\normb{v}^2}{\normltwo{v}^2}.
\end{align}

\begin{Remark}\label{WeakFormContinuity}
    Since the $\mathcal{B}$-norm is equivalent to the $H^1(\O)$-norm, the functional 
    \[
        u \mapsto \inprodb{u}{\varphi}
    \]
     belongs to the dual of $H^1(\O)$ for any fixed $\varphi \in H^1(\O)$. As a result, if $u_k \rightharpoonup u$ weakly in $H^1(\O)$ then 
    \[
        \lim_{k \to \infty} \inprodb{u_k}{\varphi} = \inprodb{u}{\varphi}
    \]
    for any $\varphi \in H^1(\O)$.
\end{Remark}

The iteration \eqref{LaplacianIterationinMainResults} in Theorem \ref{RobinConvergence} relies on solvability of the Poisson problem
\begin{equation}\label{RobinPoissonProblem}
    \begin{cases}
        -\Delta v = f & \quad \text{in } \O, \\
        v + h\partial_\nu v = 0 & \quad \text{on } \partial \O.
    \end{cases}
\end{equation}
By \cite[Theorem 6.10.10]{medkova}, if $f \in C^1(\O) \cap L^\infty(\O)$, then \eqref{RobinPoissonProblem} has a unique solution $v \in C^2(\O) \cap C^1(\overline{\O})$. Moreover, if $f \ge 0$ and $f \not\equiv 0$, then the maximum principle implies that $v$ attains its infimum at some point $x_0 \in \partial\O$. By Hopf's Lemma \cite[Section 6.4]{evans}, $\partial_\nu v(x_0) < 0$. Since $h(x_0) > 0$, we must have $v(x_0) > 0$. Thus $v$ must be strictly positive in $\overline{\O}$. Applying this observation to the iteration \eqref{LaplacianIterationinMainResults}, we see that since $u_0$ is strictly positive on $\overline{\O}$, $u_k$ is strictly positive on $\overline{\O}$ for all $k \in \N$.

\subsection{Basic Properties of the Sequence $u_k$}\label{RobinPropertiesSec}
Before we can prove Theorem \ref{RobinConvergence}, we will need some useful properties of the iterates $\{u_k\}_{k = 0}^{\infty}$ generated by \eqref{LaplacianIterationinMainResults}. When proving these properties, we will refer to the weak formulation of \eqref{LaplacianIterationinMainResults}, which is
\begin{equation}\label{LaplacianIterationWeakFormulation}
    \inprodb{u_{k + 1}}{\varphi} = R(u_k)\inprodltwo[\O]{u_k}{\varphi}, \quad \forall\varphi \in H^1(\O).
\end{equation}

\begin{Proposition}\label{Monotonicity Proposition}
    For each of the iterates $u_k$ generated by \eqref{LaplacianIterationinMainResults}, we have
    % \begin{equation}\label{monotonicityofRQandL2norm}
    %     R(u_{k+1}) ||u_{k+1}||_{L^2(\O)} \leq R(u_k) ||u_k||_{L^2(\O)},
    % \end{equation}
    % \begin{equation}\label{monotonicityofgradient}
    %     ||u_{k + 1}||_{\B} \ge ||u_k||_{\B},
    % \end{equation}
    % \begin{equation}\label{L2monotonicity}
    %     ||u_{k+1}||_{L^2(\O)} \geq ||u_{k}||_{L^2(\O)},
    % \end{equation}
    % \begin{equation}\label{monotonicityofRQ}
    %     R(u_{k+1}) \leq R(u_k).
    % \end{equation}
    \begin{align}
        R(u_{k+1}) \normltwo{u_{k+1}} &\leq R(u_k) \normltwo{u_k}, \label{monotonicityofRQandL2norm} \\
        \normb{u_{k + 1}} &\ge \normb{u_k}, \label{monotonicityofgradient} \\
        \normltwo{u_{k+1}} &\geq \normltwo{u_{k}}, \label{L2monotonicity} \\
        R(u_{k+1}) &\leq R(u_k). \label{monotonicityofRQ}
    \end{align}
\end{Proposition}

\begin{proof}[of \eqref{monotonicityofRQandL2norm}] 
    Fix $k \geq 0$. Using $u_{k+1}$ as a test function in \eqref{LaplacianIterationWeakFormulation} yields
    \[
        \normb{u_{k+1}}^2 = R(u_k) \inprodltwo{u_k}{u_{k+1}}.
    \]
    Applying the Cauchy-Schwarz inequality to the right-hand-side gives us
    \[
        \normb{u_{k+1}}^2 \leq R(u_k) \normltwo{u_k}\normltwo{u_{k+1}}.
    \]
    Using the expression for the Rayleigh quotient in \eqref{RobinRayleighQuotientShorthand}, we conclude that
    \[
        R(u_{k+1}) \normltwo{u_{k+1}} \leq R(u_k) \normltwo{u_k}.
    \]
\end{proof}

\begin{proof}[of \eqref{monotonicityofgradient}]
    Fix $k \geq 0$. Using $u_{k}$ as a test function in \eqref{LaplacianIterationWeakFormulation} yields
    \[
        \inprodb{u_{k+1}}{u_k} = R(u_k) \normltwo{u_k}^2 = \normb{u_k}^2.
    \]
    Applying the Cauchy-Schwarz inequality to the inner product $\mathcal{B}[\cdot,\cdot]$, we find that
    \[
        \normb{u_k}^2 = \inprodb{u_{k+1}}{u_k}  \leq \normb{u_k}\normb{u_{k+1}}
    \]
    Consequently,
    \[
        \normb{u_k} \leq \normb{u_{k+1}}.
    \]
\end{proof}

\begin{proof}[of \eqref{L2monotonicity}] 
     As we discussed at the beginning of this subsection, $u_k \not\equiv 0$ for any $k \in \mathbb{N}$. Thus, by \eqref{monotonicityofRQandL2norm} and \eqref{RobinRayleighQuotientShorthand}, we have
\begin{equation}\label{AnotherMonotoneQuantity}
        \frac{\normb{u_{k+1}}^2}{\normltwo{u_{k+1}}} \leq \frac{\normb{u_{k}}^2}{\normltwo{u_{k}}} \quad \text{for all } k \geq 0.
\end{equation}
   Now apply \eqref{monotonicityofgradient} to \eqref{AnotherMonotoneQuantity} to get \eqref{L2monotonicity}. 
\end{proof}

\begin{proof}[of \eqref{monotonicityofRQ}]
This follows from \eqref{L2monotonicity} and \eqref{monotonicityofRQandL2norm}.
\end{proof}

Using the monotonicity properties derived above, we can obtain an $H^1(\O)$ bound on the sequence of iterates. This will provide the compactness necessary in the proof of convergence.
\begin{Proposition}\label{uniformSobolevbounds}
    The sequence $\{ u_k \}_{k =1}^{\infty}$ is uniformly bounded in $H^1(\O)$.
\end{Proposition}

\begin{proof}
    By \eqref{monotonicityofRQandL2norm}, we have for any $k \geq 0$
    \[
        R(u_{k+1}) \normltwo{u_{k+1}} \leq R(u_0) \normltwo{u_0}.
    \]
    Since $R(u_{k+1}) \geq \lambda_R^1$, we conclude that
    \begin{equation}\label{L2bound}
        \normltwo{u_{k+1}} \leq \frac{R(u_0)}{\lambda_R^1} \normltwo{u_0}.
    \end{equation}
    Combining this with \eqref{monotonicityofRQ}, we find 
    \begin{align*}
        \normltwo{\nabla u_{k+1}}^2 \le \normb{u_{k+1}}^2 &= R(u_{k+1}) \normltwo{u_{k+1}}^2 \\
        &\leq \frac{R(u_{0})^3}{(\lambda_R^1)^{2}} \normltwo{u_{0}}^2.
    \end{align*}
 
\end{proof}

By \eqref{monotonicityofRQ}, we see that $\rho:= \lim\limits_{k \to \infty} R(u_k)$ exists and $\rho \geq \l_R^1$. It is natural to ask if $\rho = \l_R^1$; we prove this fact now.
\begin{Proposition}\label{RayleighQuotientConvergence}
    \[\lim_{k \to \infty} R(u_k) = \l_R^1.\]
\end{Proposition}
\begin{proof}
    Let $w$ be the non-negative solution of 
    \[
        \begin{cases}
            -\Delta w = \l_R^1 w & \quad \text{in } \O, \\
            w + h\partial_\nu w= 0 & \quad \text{on } \partial \O, \\
            \normltwo{w} = 1.
        \end{cases}
    \]

    Using $w$ as a test function in \eqref{LaplacianIterationWeakFormulation} yields
    \[
        \inprodb{u_{k+1}}{w} = R(u_k) \inprodltwo{u_k}{w}.
    \]
    On the other hand, we know from \eqref{RobinEigenvalueWeakFormulationShorthand} that
    \[
        \inprodb{w}{u_{k+1}} = \lambda_R^1 \inprodltwo{w}{u_{k+1}}.
    \]
    Therefore, by symmetry of the bilinear form $\mathcal{B}$,
\begin{equation}\label{eigenfunctionastestfunction}
        \lambda_R^1 \inprodltwo{u_{k+1}}{w} = R(u_k) \inprodltwo{u_k}{w}.
    \end{equation}
    Let $\alpha_k := \inprodltwo{u_k}{w}$. Since $u_k$ is strictly positive on $\overline{\O}$ for all $k \in \mathbb{N}$ and $w$ is nonnegative and not identically 0, $\alpha_k > 0$ for all $k \in \mathbb{N}$. Since $R(u_k) \geq \lambda_R^1$, we see that $\{\alpha_k\}_{k=1}^\infty$ is a non-decreasing sequence of positive real numbers. By the Cauchy-Schwarz inequality and \eqref{L2bound}, we also know that 
    \[
        \alpha_k \leq \normltwo{u_{k+1}} \leq  \frac{R(u_0)}{\lambda_R^1} \normltwo{u_0}.
    \]
    Consequently, $\alpha := \lim\limits_{k \to \infty} \alpha_k$ exists and is strictly positive. Taking the limit as $k \to \infty$ in \eqref{eigenfunctionastestfunction}, we conclude that $\lim\limits_{k \to \infty}R(u_k) = \lambda_R^1$.
\end{proof}

\setcounter{proofpart}{0}

\subsection{Convergence of the Iteration}\label{RobinProofSec}
All that is left to prove in Theorem \ref{RobinConvergence} is the convergence of $\{u_k\}_{k = 0}^{\infty}$ to a principal eigenfunction. For readability, we break the remainder of the proof of Theorem \ref{RobinConvergence} into three parts:
\begin{enumerate}
    \item First we will show that there exists a subsequence of $\{u_k\}_{k=1}^\infty$ that converges strongly in $L^2 (\O)$ and weakly in $H^1 (\O)$ to a principal eigenfunction $u_\infty$.
    \item Next we will show that any two subsequences that converge strongly in $L^2 (\O)$ and weakly in $H^1 (\O)$ will converge to the same eigenfunction.
    \item Finally, we will prove that the whole sequence converges to a principal eigenfunction strongly in $H^1(\O)$.
\end{enumerate}
\proofpart By Proposition \ref{uniformSobolevbounds}, the sequence  $\{ u_k \}_{k =1}^{\infty}$ is uniformly bounded in $H^1(\O)$. Hence, by the Rellich-Kondrachov compactness theorem, we may assume there exists a subsequence $\{ u_{k(j)} \}_{j=1}^{\infty}$ such that $u_{k(j)}$ converges weakly in $H^1(\O)$ and strongly in $L^2(\O)$ to a non-negative function $u_{\infty} \in H^1(\O)$. 

We may also choose the subsequence so that $\{u_{k(j) + 1}\}_{j=1}^\infty$ converges weakly in $H^1(\O)$ and strongly in $L^2(\O)$ to some non-negative function $w_{\infty} \in H^1(\O)$ (possibly different from $u_{\infty}$). Also note that both $u_{\infty}$ and $w_{\infty}$ are not identically zero by \eqref{L2monotonicity}.

Observe that, by \eqref{monotonicityofRQandL2norm}, we have
\begin{align*}
    R(u_{k(j+1)}) \normltwo{u_{k(j+1)}} &\leq R(u_{k(j)+1}) \normltwo{u_{k(j)+1}} \\
    &\leq  R(u_{k(j)}) \normltwo{u_{k(j)}} \quad \text{for all } j \geq 0.
\end{align*}
Since $\lim\limits_{k \to \infty} R(u_k) = \lambda_R^1$, we can take $j \to \infty$ and obtain
\begin{equation}\label{samenorms}
    \normltwo{w_{\infty}} = \normltwo{u_{\infty}}.
\end{equation}
From the weak formulation \eqref{LaplacianIterationWeakFormulation}, we have
\begin{align}\label{SubsequenceWeakFormulation}
    \inprodb{u_{k(j) + 1}}{\varphi} = R(u_{k(j)})\inprodltwo[\O]{u_{k(j)}}{\varphi} \quad \forall\varphi \in H^1(\O).
\end{align}
By Remark \ref{WeakFormContinuity}, we may take the limit as $j \to \infty$ in the weak formulation \eqref{SubsequenceWeakFormulation} using the weak convergence in $H^1(\O)$ to obtain
\begin{align}\label{uInfinityandwInfinityIndentity}
    \inprodb{w_\infty}{\varphi} = \lambda_R^1\inprodltwo{u_\infty}{\varphi} \quad \forall\varphi \in H^1(\O).
\end{align}
Choosing $\varphi = w_{\infty}$ in \eqref{uInfinityandwInfinityIndentity} and using the Cauchy-Schwarz inequality and \eqref{samenorms}, we find that
\begin{align*}
    \lambda_R^1 \normltwo{w_{\infty}}^2 & \leq R(w_{\infty}) \normltwo{w_{\infty}}^2 \\
    & =  \normb{w_\infty}^2 \\
    & = \lambda_R^1 \inprodltwo{u_{\infty}}{w_{\infty}}  \\
    & \leq \lambda_R^1  \normltwo{u_{\infty}} \normltwo{w_{\infty}} \\
    & = \lambda_R^1  \normltwo{w_{\infty}}^2.
\end{align*}

Consequently, we must have equality in the Cauchy-Schwarz inequality, and so $u_{\infty} = cw_{\infty}$ for some constant $c  > 0$ (recall that both $u_{\infty}$ and $w_{\infty}$ are non-negative). Since $\normltwo{w_{\infty}} = \normltwo{u_{\infty}}$, we conclude that $c=1$ and $u_{\infty} = w_{\infty}$. Therefore, by \eqref{uInfinityandwInfinityIndentity}, $u_{\infty}$ is a weak solution of the eigenvalue problem
\[
    \begin{cases}
        -\Delta u_{\infty} = \lambda_R^1 u_{\infty} & \quad \text{in } \O, \\
        u_{\infty} + h\partial_\nu u_{\infty}= 0 & \quad \text{on } \partial \O.
    \end{cases}
\]

\proofpart
We now want to show that any two subsequences of $\{u_k\}_{k=1}^\infty$ that converge strongly in $L^2(\O)$ and weakly in $H^1(\O)$ converge to the same limit $u_{\infty}$ obtained in the previous part.

Suppose $\{u_{k_1(j)}\}_{j=1}^\infty$ and $\{u_{k_2(j)}\}_{j=1}^\infty$ are two subsequences converging strongly in $L^2 (\O)$ and weakly in $H^1(\O)$ to $u_{1, \infty}$ and $u_{2, \infty}$, respectively. Note that Proposition \ref{RayleighQuotientConvergence} implies $R(u_{1, \infty}) = R(u_{2, \infty}) = \l_R^1$. Arguments as in Part 1 may be used to show that both $u_{1, \infty}$ and $u_{2, \infty}$ are non-negative principal eigenfunctions. 

We claim $u_{1, \infty} = u_{2, \infty}$. To do this, we construct two new subsequences $\{u_{j_1(\ell)}\}_{\ell=1}^\infty \subset \{u_{k_1(j)}\}_{j=1}^\infty$ and $\{u_{j_2(\ell)}\}_{\ell=1}^\infty \subset \{u_{k_2(j)}\}_{j=1}^\infty$ by setting $j_1(1)=k_1(1)$, then inductively choosing
\begin{align*}
    j_2(\ell) &= \min_{n \ge 1} \{k_2(n) : k_2(n) > j_1(\ell)\}, \quad &\ell \geq 1,\\
    j_1(\ell) &= \min_{n \ge 1} \{k_1(n) : k_1(n) > j_2(\ell - 1)\}, \quad &\ell \geq 2.
\end{align*}
Clearly $\{u_{j_1(\ell)}\}_{\ell=1}^\infty$ and $\{u_{j_2(\ell)}\}_{\ell=1}^\infty$ converge strongly in $L^2(\O)$ and weakly in $H^1(\O)$  to the original limits $u_{1, \infty}$ and $u_{2, \infty}$ respectively, while $j_1(\ell) < j_2(\ell)$ and $j_2(\ell)< j_1(\ell+1)$ for all $\ell$. Thus, by repeated application of the monotonicity relation \eqref{monotonicityofRQandL2norm}, we find 
\begin{align*}
    R(u_{j_2(\ell)}) \normltwo{u_{j_2(\ell)}} &\leq R(u_{j_1(\ell)}) \normltwo{u_{j_1(\ell)}}\\
    R(u_{j_1(\ell+1)}) \normltwo{u_{j_1(\ell+1)}} &\leq R(u_{j_2(\ell)}) \normltwo{u_{j_2(\ell)}}.
\end{align*}
Taking $\ell\to\infty$ in both inequalities above and then dividing by $\l_1$ yields $\normltwo{u_{1, \infty}}=\normltwo{u_{2, \infty}}$. Since both $u_{1, \infty}$ and $u_{2, \infty}$ are first eigenfunctions, they must be multiples of each other; this shows they are equal because they are both nonnegative. 

\proofpart
Note that Part 1 and Part 2 show that every subsequence that converges strongly in $L^2 (\O)$ and weakly in $H^1 (\O)$ converges to the same principal eigenfunction $u_{\infty}$ obtained in Part 1. 
By the Rellich-Kondrachov compactness theorem, this implies the full sequence $\{u_k\}_{k=1}^\infty$ converges to the principal eigenfunction $u_\infty$ strongly in $L^2 (\O)$ and weakly in $H^1 (\O)$. 

We will complete the proof of Theorem \ref{RobinConvergence} by showing that this convergence is in fact strong in $H^1(\O)$. If we subtract the weak formulations \eqref{RobinEigenvalueWeakFormulationShorthand} and \eqref{LaplacianIterationWeakFormulation} for the eigenvalue problem and the iteration, we obtain
\begin{equation*}
    \inprodb{u_{k + 1} - u_{\infty}}{\varphi} = 
    R(u_k)\inprod{u_k}{\varphi}_{L^2(\O)} - \lambda_1\inprod{u_\infty}{\varphi}_{L^2(\O)}.
\end{equation*}
Let $\varphi = u_{k + 1} - u_{\infty}$. Since $u_{k + 1} \to u_{\infty}$ in $L^2(\O)$, taking the limit as $k \to \infty$ yields
\begin{align*}
    \lim\limits_{k \to \infty}\normb{u_{k + 1} - u_{\infty}}^2 = 0.
\end{align*}
By the equivalence of $\normb{\cdot}$ and $\normhone[\O]{\cdot}$, $u_k$ converges to $u_{\infty}$ strongly in $H^1(\O)$.%
\hspace*{\fill}
$\Box$

\section{Inverse Iteration for the Mixed Eigenvalue Problem}\label{sec:Mixed}

\subsection{Background}\label{MixedBackgroundSec} Again, by adapting the techniques in \cite[Section 6.5, Theorems 1-2]{evans} and \cite[Theorems 7.41.5-7.41.6]{medkova}, we have the following:
\begin{enumerate}%[label = (\roman*)]
    \item The eigenvalues in \eqref{MixedEigenvalueProblem} are discrete, positive, and diverge to infinity.
    \item The principal eigenvalue admits the variational characterization
    \[\lambda_M^1 = \min_{v\in \mathcal{C}} \r(v),\]
    where the Rayleigh quotient $\r(v)$ is as in \eqref{MixedRayleighQuotient} and $\mathcal{C}$ is as in \eqref{MixedTestFunctionClass}. 
    \item The principal eigenfunction has a sign on $\O$.
    \item The principal eigenvalue $\lambda^1_M$ is simple.
\end{enumerate}

We say that $u \in \mathcal{C}$ is a weak solution of \eqref{MixedEigenvalueProblem} with $\lambda = \lambda_M^1$ if it satisfies
\begin{equation}\label{MixedWeakFormulation}
    \int_{\Omega} \nabla u \cdot \nabla \varphi \, dx = \l_M^1 \int_{\Omega} u \varphi \, dx \qquad \forall \varphi \in \mathcal{C}.
\end{equation}
We will first show that $\mathcal{C}$ is closed under weak convergence.
\begin{Lemma}\label{CWeaklyClosed}
    If $\{u_k\}_{k=1}^\infty \subseteq \mathcal{C}$ converges weakly to $u \in H^1 (\O)$, then $u \in \mathcal{C}$.
\end{Lemma}

\begin{proof}
     Let $\{u_k\}_{k=1}^\infty \subseteq \mathcal{C}$ be such that $u_k$ converges weakly to $u \in H^1(\O)$. For fixed $\varphi \in L^2(\GD)$, we let $L_\varphi$ be the linear functional on $H^1(\O)$ defined as
     \begin{align*}
         L_\varphi(v) := \int_{\GD}v\varphi \, d\sigma.
     \end{align*}
     Note that $L_\varphi$ is bounded for all $\varphi \in L^2(\GD)$ since
     \begin{align*}
         |L_\varphi(v)| &= \left|\int_{\GD}v\varphi d\sigma\right| \\
         &\le \norm{v}_{L^2(\GD)}\norm{\varphi}_{L^2(\GD)} \\
         &\le \norm{v}_{L^2(\partial\O)}\norm{\varphi}_{L^2(\GD)} \\
         &\le C \norm{v}_{H^1(\O)}\norm{\varphi}_{L^2(\GD)},
     \end{align*}
     where the last inequality follows from the trace theorem \cite[Section 5.5]{evans}. Therefore
     \begin{align*}
         \lim\limits_{k \to \infty}L_\varphi(u_k) = L_\varphi(u).
     \end{align*}
     But $L_\varphi(u_k) = 0$ for all $k \in \mathbb{N}$ because $u_k \in \mathcal{C}$. Thus, $L_\varphi(u) = 0$. Since this holds for all $\varphi \in L^2(\GD)$ we must have $u = 0$ almost everywhere on $\GD$, i.e. $u \in \mathcal{C}$.
\end{proof}

The iteration \eqref{MixedIteration} in Theorem \ref{MixedConvergence}, relies on the solvability of the Poisson problem with mixed boundary conditions
\begin{equation}\label{PoissonMixed}
\begin{cases}
    -\Delta v = f &\quad \text{in } \O,\\
    v = 0 &\quad \text{on } \GD, \\
    \partial_\nu v = 0 &\quad \text{on } \GN.
\end{cases}
\end{equation}
By \cite[Theorem 6.20.5]{medkova}, if $f \in C^1(\O) \cap L^\infty(\O)$ and $\O$ is bounded and smooth then \eqref{PoissonMixed} has a unique classical solution $v \in C^2(\O) \cap C^1(\overline \O)$. Then by the strong maximum principle, either $v$ attains an infimum of $0$ on $\Gamma_D$ or it attains its infimum at some $x_0 \in \GN$. In the second case, $\partial_\nu v(x_0) < 0$, which is a contradiction, so the minimum must be attained on $\Gamma_D$ and so $v \ge 0$ in $\overline{\O}$. Applying this observation to the iteration \eqref{MixedIteration}, we see that since $u_0$ is strictly positive on $\O$, $u_k$ is non-negative in $\overline{\O}$ for all $k \in \N$.

\subsection{Basic Properties of the Sequence $u_k$}\label{MixedPropertiesSec}
Before we can prove Theorem \ref{MixedConvergence}\, we will need some useful properties of the iterates $\{u_k\}_{k=0}^\infty$ generated by \eqref{MixedIteration}. When proving these properties, we will refer to the weak formulation of \eqref{MixedIteration}, which is 
\begin{equation}\label{MixedIterationWeakFormulation}
\int_{\O} \nabla u_{k+1} \cdot \nabla\varphi \, dx = \r(u_k) \int_\O u\varphi \, dx \qquad \forall \varphi \in \mathcal{C}.
\end{equation}

\begin{Proposition}\label{MonotonicityPropositionMixed}
    For each of the iterates $u_k$ generated by \eqref{MixedIteration}, we have
    % \begin{equation}\label{monotonicityofRQandL2normMixed}
    %     \r(u_{k+1}) ||u_{k+1}||_{L^2(\O)} \leq \r(u_k) ||u_k||_{L^2(\O)}.
    % \end{equation}
    % \begin{equation}\label{monotonicityofgradientMixed}
    %     \normltwo{\nabla u_{k+1}}  \geq \normltwo{\nabla u_{k}}.
    % \end{equation}
    % \begin{equation}\label{L2monotonicityMixed}
    %     ||u_{k+1}||_{L^2(\O)} \geq ||u_{k}||_{L^2(\O)}.
    % \end{equation}
    % \begin{equation}\label{monotonicityofRQMixed}
    %     \r(u_{k+1}) \leq \r(u_k).
    % \end{equation}
    \begin{align}
        \r(u_{k+1}) \normltwo{u_{k+1}} &\leq \r(u_k) \normltwo{u_k}, \label{monotonicityofRQandL2normMixed} \\
        \normltwo{\nabla u_{k+1}}  &\geq \normltwo{\nabla u_{k}}, \label{monotonicityofgradientMixed} \\
        \normltwo{u_{k+1}} &\geq \normltwo{u_{k}}, \label{L2monotonicityMixed} \\
        \r(u_{k+1}) &\leq \r(u_k). \label{monotonicityofRQMixed}
    \end{align}
\end{Proposition}

The proofs below are very similar to those in Section \ref{RobinPropertiesSec}, but the bilinear form in the weak formulation \eqref{MixedIterationWeakFormulation} is not an inner product, so we have elected to provide the necessary details. \\
\begin{proof}[of \eqref{monotonicityofRQandL2normMixed}] 
    Fix $k \geq 0$. Using $u_{k+1}$ as a test function in \eqref{MixedIterationWeakFormulation} yields
    \[
        \int_\O |\nabla u_{k+1}|^2\, dx = \r(u_k) \int_\O u_k u_{k+1}\, dx.
    \]
    Applying the Cauchy-Schwarz inequality to the right-hand-side gives us
    \[
        \int_\O |\nabla u_{k+1}|^2\, dx \leq \r(u_k) \normltwo{u_k}\normltwo{u_{k+1}}.
    \]
    Using the definition of $\r(u_{k+1})$, we conclude that
    \[
        \r(u_{k+1}) \normltwo{u_{k+1}} \leq \r(u_k) \normltwo{u_k}.
    \]
\end{proof}

\begin{proof}[of \eqref{monotonicityofgradientMixed}]
    Fix $k \geq 0$. Using $u_{k}$ as a test function in \eqref{MixedIterationWeakFormulation} yields
    \[
        \int_\O \nabla u_{k+1} \cdot \nabla u_k \, dx = \r(u_k) \normltwo{u_k}^2.
    \]
    Applying the Cauchy-Schwarz inequality to the left-hand-side yields
    \[
        \normltwo{\nabla u_{k}}\normltwo{\nabla u_{k+1}} \geq  \r(u_k) \normltwo{u_k}^2 = \normltwo{\nabla u_{k}}^2.
    \]
    Consequently,
    \[
        \normltwo{\nabla u_{k}} \leq \normltwo{\nabla u_{k+1}}.
    \]
\end{proof}

\begin{proof}[of \eqref{L2monotonicityMixed}] 
    We first observe that by \eqref{monotonicityofRQandL2normMixed} and the definition of the Rayleigh quotient, we have
    \begin{equation}\label{AnotherMonotoneQuantityMixed}
        \frac{\normltwo{\nabla u_{k+1}}^2}{\normltwo{u_{k+1}}} \leq \frac{\normltwo{\nabla u_{k}}^2}{\normltwo{u_{k}}} \quad \text{for all } k \geq 0.
    \end{equation}
    We may note here that $u_k \not\equiv 0$ as long as $u_0 \not\equiv 0$ as we have shown that $\{\normltwo{\nabla u_k}^2\}_{k=1}^\infty$ is an increasing sequence.
    Now apply \eqref{monotonicityofgradientMixed} to \eqref{AnotherMonotoneQuantityMixed} to get \eqref{L2monotonicityMixed}. 
\end{proof}

\begin{proof}[of \eqref{monotonicityofRQMixed}]
This follows from \eqref{L2monotonicityMixed} and \eqref{monotonicityofRQandL2normMixed}.
\end{proof}

Using the monotonicity properties above, we can derive an $H^1$ bound on the sequence of iterates.
\begin{Proposition}\label{uniformSobolevboundsMixed}
    The sequence $\{ u_k \}_{k =1}^{\infty}$ is uniformly bounded in $H^1(\O)$.
\end{Proposition}

\begin{proof}
Using the monotonicity properties established in Proposition \ref{MonotonicityPropositionMixed}, an identical argument to the proof of Proposition \ref{uniformSobolevbounds} shows that
    \begin{equation}\label{L2boundMixed}
         \normltwo{u_{k+1}} \leq \frac{\r(u_0)}{\lambda_M^1} \normltwo{u_0},
     \end{equation}
     and
    \begin{align*}
        \normltwo{\nabla u_{k+1}}^2 & \leq \frac{\r(u_{0})^3}{(\lambda_M^1)^{2}}\normltwo{u_{0}}^2.
    \end{align*}
\end{proof}

By \eqref{monotonicityofRQMixed}, we see that $\rho := \lim\limits_{k \to \infty} \r(u_k)$ exists and $\rho \geq \lambda_M^1$. It is natural to ask if $\rho = \lambda_M^1$; we show this to be true and the proof is quite similar to that of Proposition \ref{RayleighQuotientConvergence} in Section \ref{RobinPropertiesSec}.

\begin{Proposition}\label{RayleighQuotientConvergenceMixed}
    \[\lim_{k \to \infty} \r(u_k) = \l_M^1.\]
\end{Proposition}
\begin{proof}
    Let $w \in \mathcal{C}$ be the non-negative solution of 
    \[
        \begin{cases}
            -\Delta w = \l_M^1 w & \quad \text{in } \O, \\
            w = 0 & \quad \text{on } \GD, \\
            \partial_\nu w = 0 &\quad \text{on } \GN,\\
            \normltwo{w} = 1.
        \end{cases}
    \]
    In particular, $w$ is strictly positive in $\O$ in accordance with the spectral theory outlined in Section \ref{MixedBackgroundSec}. Using $w$ as a test function in \eqref{MixedIterationWeakFormulation} yields
    \[
        \int_\O \nabla w \cdot \nabla u_{k+1} \, dx = \r (u_k) \inprodltwo{u_k}{w}.
    \]
    On the other hand, we know from the weak formulation \eqref{MixedWeakFormulation} for $w$ that
    \[
        \int_\O \nabla w \cdot \nabla u_{k+1} \, dx = \lambda_M^1 \inprodltwo{w}{u_{k+1}}.
    \]
    Therefore,
    \begin{equation}\label{eigenfunctionastestfunctionMixed}
        \lambda_M^1 \inprodltwo{u_{k+1}}{w} = \r (u_k) \inprodltwo{u_k}{w}.
    \end{equation}
    Let $\alpha_k := \inprodltwo{u_k}{w}$. Since $u_k$ and $w$ are both strictly positive on $\O$, $\alpha_k > 0$ for all $k$. Since $R(u_k) \geq \lambda_R^1$, we see that $\{\alpha_k\}_{k=1}^\infty$ is an non-decreasing sequence of positive real numbers. By the Cauchy-Schwarz inequality and \eqref{L2boundMixed}, we also know that 
    \[
        \alpha_k \leq \normltwo{u_{k+1}} \leq \frac{\r (u_0)}{\lambda_M^1} \normltwo{u_0}.
    \]
    Consequently, $\alpha := \lim\limits_{k \to \infty} \alpha_k$ exists and is strictly positive, so taking the limit as $k \to \infty$ in \eqref{eigenfunctionastestfunctionMixed} we conclude that $\lim\limits_{k \to \infty}\r (u_k) = \lambda_M^1$.
\end{proof}

\subsection{Convergence of the Iteration}\label{MixedProofSec}
We break up the remainder of the proof into three parts, as in Section \ref{RobinProofSec}:
\begin{enumerate}
    \item First we will show that there exists a subsequence of $\{u_k\}_{k=1}^\infty$ that converges strongly in $L^2 (\O)$ and weakly in $H^1 (\O)$ to a principal eigenfunction $u_\infty$.
    \item Next we will show that any two subsequences that converge strongly in $L^2 (\O)$ and weakly in $H^1 (\O)$ will converge to the same eigenfunction.
    \item Finally, we will prove that the whole sequence converges to a principal eigenfunction strongly in $H^1(\O)$
\end{enumerate}
With the properties proven in Section \ref{MixedPropertiesSec}, the proofs of Parts 1 and 2 follow nearly identically to that of Theorem \ref{RobinConvergence} detailed in Section \ref{RobinProofSec}. The same test functions are all valid because $\mathcal{C}$ is weakly closed by Lemma \ref{CWeaklyClosed}. However, the proof of Part 3 is different enough that we provide some details below.

\begin{proof}[of Part 3]
    Let $u_\infty$ denote the eigenfunction to which the sequence of iterates $\{u_k\}_{k=1}^\infty$ converges to strongly in $L^2 (\O)$ and weakly in $H^1 (\O)$. If we subtract the weak formulations \eqref{MixedWeakFormulation} and \eqref{MixedIterationWeakFormulation} for the eigenvalue problem and the iteration, we obtain
    \begin{equation*}\label{DifferenceWeakFormulation}
        \int_\O \nabla (u_{k+1} - u_\infty) \cdot \nabla \varphi \, dx = 
        \r (u_k)\int_{\O}u_k\varphi \, dx - \lambda_M^1\int_{\O}u_\infty\varphi \, dx.
    \end{equation*}
    Plugging in $\varphi = u_{k + 1} - u_\infty \in \mathcal{C}$ and taking the limit as $k \to \infty$ on both sides yields
    \begin{align*}
        \lim\limits_{k \to \infty}\normltwo{\nabla(u_{k + 1} - u_\infty)} = 0.
    \end{align*}
    Since $\lim\limits_{k \to \infty}\normltwo{u_{k + 1} - u_\infty} = 0$ by strong $L^2(\O)$ convergence, the desired claim follows.
\end{proof}

\section{Inverse Iteration for an Eigenvalue Problem in Optimal Insulation}\label{sec:Double}

\subsection{Background}

We now propose and investigate a promising inverse iterative scheme for a problem in optimal insulation studied by Bucur, Buttazo, and Nitsch in \cite{BucurButtazzo}. For additional background on this optimization problem, we refer to \cite{TwoOptProb} and \cite[Section 7.2.3]{HenrotBook}, and to \cite{BartelsButtazzo} for numerical results.

Recall the Rayleigh quotient \eqref{RobinRayleighQuotient} for the Robin eigenvalue problem:
\[
    R(v, h) = \frac{\int_{\O} |\nabla v|^2 \, dx + \int_{\partial \O} h^{-1} v^2 \, d\sigma}{\int_{\O} v^2 \, dx}.
\]
The dependency on $h$ will be kept explicit for the purposes of this discussion. Given a non-negative function $h$ on $\partial \O$, denote the principal Robin eigenvalue by $\lambda_R^1(h)$. The optimization problem studied by \cite{BucurButtazzo} is to minimize  $\lambda_R^1(h)$ among all $h \in \mathcal{H}_m$, where
\[
    \mathcal{H}_m := \left\{h : \partial \O \rightarrow [0, \infty) : \int_{\partial \O} h \, d\sigma = m \right\}
\]
is the class of insulation profiles of fixed mass $m> 0$. Let $\lambda_m := \min\limits_{h \in \mathcal{H}_m} \lambda_R^1(h)$. If $u_h$ is a principal eigenfunction for the Robin eigenvalue problem \eqref{RobinEigenvalueProblem}, we can write
\begin{align*}
    \lambda_m & = \min_{h \in \mathcal{H}_m} \lambda_R^1(h)  \\
    & = \min_{h \in \mathcal{H}_m} R(u_h, h) \\
    & = \min_{h \in \mathcal{H}_m} \; \min_{v \in H^1(\O)}  R(v, h) \\
    & = \min_{v \in H^1(\O)}\;   \min_{h \in \mathcal{H}_m}  R(v, h).
\end{align*}
Now, for any $v \in H^1(\O)$ fixed, we have
\begin{align*}
    \min_{h \in \mathcal{H}_m}   R(v, h) & = \min_{h \in \mathcal{H}_m}   \frac{ \int_{\O} |\nabla v|^2 + \int_{\partial \O} h^{-1} v^2 \ d\sigma}{\int_{\O} v^2 \ dx} \\
    & =\frac{ \int_{\O} |\nabla v|^2 + \min_{h \in \mathcal{H}_m} \left\{\int_{\partial \O} h^{-1} v^2 \ d\sigma \right\}}{\int_{\O} v^2 \ dx}
\end{align*}
By the Cauchy-Schwarz inequality, we have for any positive $h$ and $v \in H^1(\O)$
\begin{equation*}
    \left(\int_{\partial \O} h  \ d\sigma \right) \left( \int_{\partial \O} \frac{v^2}{h}  \ d\sigma \right) \geq \left(\int_{\partial \O} |v|  \ d\sigma \right)^2.
\end{equation*}
Since $h \in \mathcal{H}_m$, we have 
\begin{equation}\label{CSDI2}
     \left( \int_{\partial \O} \frac{v^2}{h}  \ d\sigma \right) \geq \frac{1}{m} \left(\int_{\partial \O} |v|  \ d\sigma \right)^2.
\end{equation}
Equality in \eqref{CSDI2} is achieved only when we have equality in the Cauchy-Schwarz inequality. This means $\int_{\partial \O} h^{-1} v^2 \ d\sigma$ achieves a minimum over $h \in \mathcal{H}_m$ if and only if
\[
    h = \frac{m |v|}{\int_{\partial \O} |v|  \ d\sigma}.
\]
Substituting this expression for $h$ into the Rayleigh quotient $R(v,h)$, we find that 
\[
    \lambda_m =  \min_{v \in H^1(\O)} R_m(v),
\]
where 
\begin{equation*}
    R_m(v)  :=\frac{ \int_{\O} |\nabla v|^2 + \frac{1}{m} \left(\int_{\partial \O} |v| \ d\sigma\right)^2}{\int_{\O} v^2 \ dx}.
\end{equation*}
By \cite[Section 3]{BucurButtazzo}, the existence of minimizers of $R_m(v)$ follows via the direct method in calculus  of variations.  Additionally, if $\hat{u}$ is a minimizer of $R_m(v)$, then $\lambda_m =  \min\limits_{h \in \mathcal{H}_m} \lambda_R^1(h)$ is attained by the insulation profile \[\hat{h} := \frac{m |\hat{u}|}{\int_{\partial \O} |\hat{u}|  \ d\sigma}.\]
Finally, by \cite[Remark 2.5]{BucurButtazzo}, any minimizer of $R_m(v)$ solves the nonlocal eigenvalue problem 
\begin{equation*}
    \begin{cases}
        -\Delta u = \lambda_m u & \quad \text{in } \O,\\
        0 \in m \partial_{\nu} u+ H(u) \int_{\partial \O} |u| \ d\sigma & \quad \text{on } \partial \O
    \end{cases}
\end{equation*}
where $H(t)$ is the multi-valued mapping
\[
    H(t) \in  
    \begin{cases}
        \{1\}& \quad \text{if } t > 0,\\
        \{-1\}& \quad \text{if } t < 0,\\
        [-1,1]& \quad \text{if } t =0.
    \end{cases}
\]

\subsection{Inverse Iteration}
We wish to develop an iterative approach to finding a minimizer $\hat{u}$ of $R_m(\cdot)$. Initialize the iteration by choosing $u_0 \in C^1(\O) \cap L^\infty(\O)$ strictly positive on $\overline{\O}$ and let 
\[h_0 := \frac{m |u_0|}{\int_{\partial \O} |u_0|  \ d\sigma}.\] 
We now consider the following iteration:
\begin{equation}\label{RobinEigenvalueIteration}
    \begin{cases}
        -\Delta u_{k+1} = R_m (u_k) u_k & \quad \text{in } \O,\\
        u_{k+1} + h_k \partial_{\nu} u_{k+1} = 0 & \quad \text{on } \partial \O. 
    \end{cases}
    \qquad
    h_{k+1} := \frac{m |u_{k+1}|}{\int_{\partial \O} |u_{k+1}|  \ d\sigma}
\end{equation}
Because $u_0$ is strictly positive, $h_0$ is strictly positive. Therefore, by the regularity theory for the Robin Poisson Problem \eqref{RobinPoissonProblem}, there exists a unique $u_1 \in C^2(\O) \cap C^1(\overline{\O})$, strictly positive on $\overline{\O}$, satisfying \eqref{RobinEigenvalueIteration}. Then $h_1$ is also strictly positive. We can repeat this argument to deduce that \eqref{RobinEigenvalueIteration} is well-defined and that $u_k$ and $h_k$ are always strictly positive.

Similar to \eqref{LaplacianIterationinMainResults} and \eqref{MixedIteration}, the iterates of \eqref{RobinEigenvalueIteration} also satisfy several monotonicity properties, whose proofs we will provide below. Before doing so, we note that the following consequence of the Cauchy-Schwarz inequality: for all $v \in H^1(\O)$ and $h \in \mathcal{H}_m$, we must have 
\[
    R_m(v) \leq R(v,h) \quad \text{ with equality if and only if } \quad h = \frac{m |v|}{\int_{\partial \O} |v| \, d\sigma}.
\]
In particular, we have
\[
    R_m(u_k) = R(u_k, h_k).
\]
We will use the weak form of \eqref{RobinEigenvalueIteration}:
\begin{equation}\label{RobinIterationWeakForm}
    \int_{\O} \nabla u_{k+1} \cdot \nabla \varphi \, dx + \int_{\partial \O} \frac{u_{k+1} \varphi}{h_k} \, d\sigma = R_m(u_k) \int_{\O} u_k \varphi\, dx \quad \forall \varphi \in H^1(\O).
\end{equation}
We can restate this in terms of the bilinear forms
\begin{equation*}
    B_k [v,\varphi] := \int_\Omega  \nabla v \cdot \nabla \varphi \ dx + \int_{\partial \Omega} \frac{v\varphi}{h_k} \; d\sigma
\end{equation*}
so that \eqref{RobinIterationWeakForm} is equivalent to 
\[
    B_k[u_{k+1},\varphi] = R_m(u_k) \left\langle u_k, \varphi \right\rangle_{L^2(\O)} \quad \forall \varphi \in H^1(\O).
\]

\begin{Proposition}
Along a sequence of iterates generated by \eqref{RobinEigenvalueIteration}, we have the following properties for all $k \geq 0$:
\begin{align}
    R_m(u_{k + 1})\normltwo{u_{k + 1}} &\le R_m(u_k)\normltwo{u_k}, \label{monotonicity1}\\
    \normltwo{u_{k}} &\leq \normltwo{u_{k+1}}, \label{monotonicity3} \\
    R_m (u_{k+1}) &\le R_m (u_k). \label{monotonicity4}
\end{align}
\end{Proposition}

\begin{proof}[of \eqref{monotonicity1}] 
    We use $\varphi = u_{k+1}$ as a test function in \eqref{RobinIterationWeakForm} to obtain
\begin{align*}
    B_k[u_{k+1},u_{k+1}] & = R_m(u_k) \left\langle u_k, u_{k+1} \right\rangle_{L^2(\O)} \\
    & \le R_m(u_k)\normltwo{u_k}\normltwo{u_{k + 1}}. 
\end{align*}
Since $B_k[u_{k+1},u_{k+1}] = R(u_{k + 1}, h_k)\normltwo{u_{k + 1}}^2$, dividing both sides through by $\normltwo{u_{k + 1}}$ yields
\begin{align*}
    R(u_{k + 1}, h_k)\normltwo{u_{k + 1}} \le R_m(u_k)\normltwo{u_k}.
\end{align*}
Then, the desired result is a consequence of the fact that $R_m (u_{k+1}) \leq R(u_{k+1}, h_k)$. 
\end{proof}

\begin{proof}[of \eqref{monotonicity3}]
    Using $\varphi = u_k$ as a test function in \eqref{RobinIterationWeakForm}, we obtain
\begin{align*}
    B_k[u_{k+1}, u_k] & = R_m(u_k) \normltwo{u_k}^2 \\
    & = R(u_k,h_k) \normltwo{u_k}^2  \\
    & = B_k[u_k,u_k].
\end{align*}
Note that, since $h_k$ is strictly positive on $\partial \O$,  $B_k[\cdot, \cdot]$ defines an inner product on $H^1(\O)$, so by the Cauchy-Schwarz inequality
\begin{align*}
    B_k [u_k, u_k] = B_k[u_{k+1}, u_k] \leq (B_k [u_k,u_k])^{1/2} (B_k [u_{k+1},u_{k+1}])^{1/2}.
\end{align*}
Due to the positivity of the iterates $u_k$, we may divide both sides by $B_k[u_k, u_k]$ then square both sides to obtain
\[
    B_k [u_k, u_k] \le B_k [u_{k+1},u_{k+1}].
\]
It is then immediate by our weak formulations that 
\[
    R(u_k,h_k) \normltwo{u_k}^2  \leq R(u_{k+1},h_k) \normltwo{u_{k+1}}^2.
\]
Now, after applying \eqref{monotonicity1}, we have
\[
    R(u_{k+1},h_k) \normltwo{u_{k+1}}^2  \leq R(u_k, h_k) \normltwo{u_k} \normltwo{u_{k+1}}.
\]
Therefore, $\normltwo{u_k} \leq \normltwo{u_{k+1}}$.
\end{proof}
\begin{proof}[of \eqref{monotonicity4}]
     This follows from \eqref{monotonicity1} and \eqref{monotonicity3}.
\end{proof}

\begin{Proposition}
    The sequence $\{u_k\}_{k=1}^\infty$ is uniformly bounded in $H^1 (\O)$.
\end{Proposition}
\begin{proof}
 By \eqref{monotonicity1} and the fact that $\lambda_m \leq R_m(v)$ for all $v \in H^1(\O)$, we have
\begin{equation}\label{monotonicity2}
    \normltwo{u_{k+1}} \leq \frac{R_m(u_0)}{\lambda_m} \normltwo{u_0}.
    \end{equation}
    By \eqref{monotonicity4}, we have for each $k \geq 1$
    \[
        \normltwo{\nabla u_k}^2 \leq R_m (u_k) \normltwo{u_k}^2 \leq R_m (u_0) \normltwo{u_k}^2.
    \]
    Then, by \eqref{monotonicity2}
    \[
        \normltwo{\nabla u_k}^2 \leq \frac{R_m(u_0)^3}{\l_m^2} \normltwo{u_0}^2.
    \]
\end{proof}

\subsection{Remarks on the Convergence of the Iteration}

The properties developed in the previous subsection suggest that the arguments in Section \ref{RobinProofSec} can be adapted to show that the iterates generated by \eqref{RobinEigenvalueIteration} converge to a minimizer of $R_m$. Unfortunately, our attempts to  prove such a convergence result were unsuccessful, primarily due to the fact that the boundary conditions change at each step of the iteration because of the update to $h_k$ in \eqref{RobinEigenvalueIteration}. A more subtle issue is the fact that minimizers of $R_m$ may not be unique. Indeed, a major result of \cite{BucurButtazzo} is the following symmetry-breaking phenomenon: if $\O$ is the unit ball, there exists $m_0 > 0$ such that for $m > m_0$, the minimizer of $R_m$ is unique and radial, while for $m < m_0$, any minimizer  of $R_m$ must be non-radial. Consequently, the optimal insulating layer is uniform if $m > m_0$ and non-uniform if $m < m_0$. The intuitive explanation for this phenomenon, as outlined in \cite[Remark 3.2]{BucurButtazzo}, is that it is better to distribute the insulation in a non-uniform manner if the total amount of insulation is small.

We conjecture that a judicious choice of initial function $u_0$ will guarantee convergence of the iterative scheme \eqref{RobinEigenvalueIteration}, at least when $R_m$ has a unique minimizer. For instance, if $u_0$ is a radial function, then the iterates $u_k$ must also be radial, and so cannot converge to a minimizer of $R_m$ when $m < m_0$. %Numerical examples with different choices of $u_0$ might be helpful in identifying the properties necessary to guarantee convergence of \eqref{RobinEigenvalueIteration}; we plan to investigate this in future work.
Numerical examples with different choices of $u_0$ might be helpful in identifying the properties necessary to guarantee convergence of \eqref{RobinEigenvalueIteration}, and this is the subject of ongoing investigation.

\section*{Acknowledgments}

This work was carried out while the authors were participating in the 2024 Mathematics Summer REU Program at Lafayette College; we acknowledge support from the NSF. We would like to thank our REU co-mentors, Profs. Farhan Abedin and Jun Kitagawa, for suggesting the problems tackled in this paper and for their guidance. We would also like to thank Profs. Jeff Liebner and Qin Lu for organizing the Lafayette REU which was an exceptional research environment.

\bibliographystyle{amsplain}
\bibliography{refs}
\end{document}